\documentclass[12pt, twocolumn]{IEEEtran}
\usepackage{eurosym}
\usepackage{amsfonts}
\usepackage{amsmath}
\usepackage[section]{placeins}
\usepackage{graphicx}
\usepackage{amsthm}
\usepackage[]{mcode}

\begin{document}
\title{Sparse Solutions of an Undetermined Linear System}
\author{Maddullah Almerdasy \\
New York University Tandon School of Engineering
}

\maketitle

\begin{abstract}
This work proposes a research problem of finding
sparse solution of undetermined Linear system with some applications.
Two approaches how to solve the compressive sensing problem: using $l_1$ approach , the $l_q$ approach
with $0<q \leq 1$. Compressive sensing algorithms are designed to cope with
ambiguities introduced by under-sampling. We propose an algorithm for restricted Isometry and how it can be used to constrain
the undetermined linear system to eventually get a unique solution.
\end{abstract}

\section{System Model and Background}

Compressive sensing (CS) aims to solve the following optimization problem to estimate the sparsest vector $\mathbf{x}$

\begin{equation}
\begin{split}
\text{min}_\mathbf{x}\;\; & \|\mathbf{x}\|_0 \\
\text{subject to}\;\; & y_i = \mathbf{b}^H_i \mathbf{x}, \; \; \mathbf{y}\in \mathbb{C}^N \; \mathbf{b}_i \in \mathbb{C}^n \; \mathbf{x} \in \mathbb{C}^n\\
& i = 1,\ldots,N \label{CS}
\end{split}
\end{equation}

where the $\|\mathbf{x}\|_0$ counts the number of nonzeros in $\mathbf{x}$. This problem is known to be NP-hard, and the theory of compressive sensing has found many approximations to solve this problem. One such approximation is known as \emph{basis pursuit} (BP) and solves the following convex approximation to (\ref{CS})

\begin{equation}
\begin{split}
\text{min}_\mathbf{x}\;\; & \|\mathbf{x}\|_1 \\
\text{subject to}\;\; & y_i = \mathbf{b}^H_i \mathbf{x},\;\; \mathbf{y}\in \mathbb{C}^N \; \mathbf{b}_i \in \mathbb{C}^n \; \mathbf{x} \in \mathbb{C}^n \\
& i = 1,\ldots,N \label{BP}.
\end{split}
\end{equation}

This approach allows one to solve for the sparsest vector $\mathbf{x}$ for a linear model under certain conditions, which have been studied extensively \cite{BP}. Recently, there has been an interest to generalize CS to nonlinear models. In particular in this report we will focus on solving

\begin{equation}
\begin{split}
\text{min}_\mathbf{x}\;\; & \|\mathbf{x}\|_0 \\
\text{subject to}\;\; & y_i = a_i + \mathbf{b}^H_i \mathbf{x} + \mathbf{x}^H \mathbf{c}_i + \mathbf{x}^H \mathbf{Q}_i \mathbf{x} \\
& i = 1,\ldots,N \label{QCS}
\end{split}
\end{equation}

where, $\mathbf{y}\in \mathbb{C}^N, \; \mathbf{a}_i \in \mathbb{C}^n, \; \mathbf{b}_i \in \mathbb{C}^n, \;  \mathbf{c}_i \in \mathbb{C}^n, \; \mathbf{Q} \in \mathbb{C}^{n\times n} \; \text{and }\mathbf{x}  \in \mathbb{R}^n.$ Solving (\ref{QCS}) will allow us to solve CS problems where the model takes a quadratic form rather than a linear form.

\section{Motivation}
In this section we give reasons why we want to solve the sparse solution of undetermined systems of linear equations.
\subsection{Signal and image Compression}
This is the most direct application. Lets consider a equation,
$y = Hx$ in which $y$ is an under sampled version of image $x$
in a fourier domain. if we no before hand that $x$ is sparse
Image, we can recover $x$ from $y$ using equation that provides
an approximation of $x$ from $y$.

\subsection{Compressed Sensing}
Compressive sensing (CS) aims to solve the following optimization problem to estimate the sparsest vector $\mathbf{x}$

\begin{equation}
\begin{split}
\text{min}_\mathbf{x}\;\; & \|\mathbf{x}\|_0 \\
\text{subject to}\;\; & y_i = \mathbf{b}^H_i \mathbf{x}, \; \; \mathbf{y}\in \mathbb{C}^N \; \mathbf{b}_i \in \mathbb{C}^n \; \mathbf{x} \in \mathbb{C}^n\\
& i = 1,\ldots,N \label{CS}
\end{split}
\end{equation}

where the $\|\mathbf{x}\|_0$ counts the number of nonzeros in $\mathbf{x}$. This problem is known to be NP-hard, and the theory of compressive sensing has found many approximations to solve this problem. One such approximation is known as \emph{basis pursuit} (BP) and solves the following convex approximation to (\ref{CS})

\begin{equation}
\begin{split}
\text{min}_\mathbf{x}\;\; & \|\mathbf{x}\|_1 \\
\text{subject to}\;\; & y_i = \mathbf{b}^H_i \mathbf{x},\;\; \mathbf{y}\in \mathbb{C}^N \; \mathbf{b}_i \in \mathbb{C}^n \; \mathbf{x} \in \mathbb{C}^n \\
& i = 1,\ldots,N \label{BP}.
\end{split}
\end{equation}

This approach allows one to solve for the sparsest vector $\mathbf{x}$ for a linear model under certain conditions, which have been studied extensively \cite{BP}. Recently, there has been an interest to generalize CS to nonlinear models. In particular in this report we will focus on solving

\begin{equation}
\begin{split}
\text{min}_\mathbf{x}\;\; & \|\mathbf{x}\|_0 \\
\text{subject to}\;\; & y_i = a_i + \mathbf{b}^H_i \mathbf{x} + \mathbf{x}^H \mathbf{c}_i + \mathbf{x}^H \mathbf{Q}_i \mathbf{x} \\
& i = 1,\ldots,N \label{QCS}
\end{split}
\end{equation}

where, $\mathbf{y}\in \mathbb{C}^N, \; \mathbf{a}_i \in \mathbb{C}^n, \; \mathbf{b}_i \in \mathbb{C}^n, \;  \mathbf{c}_i \in \mathbb{C}^n, \; \mathbf{Q} \in \mathbb{C}^{n\times n} \; \text{and }\mathbf{x}  \in \mathbb{R}^n.$ Solving (\ref{QCS}) will allow us to solve CS problems where the model takes a quadratic form rather than a linear form.

\subsection{Error Correction Code}
Let $z$ be a vector which has encoded $x$ by a linear system. A of size $m\times n$ where $m>n$. That is $z = Ax$. Now if the
the $z$ obtained is corrupted by noise such that $z = Ax + e$. In such a case we need to obtain the error corrupting $z$, lets
consider a matrix $B$ which is $n \times m$ such that $BA = 0$ hence we obtain $Bz = Be$, we consider $y = Bz$ therefore we get a linear equation $y = Be$ we can use (1) to solve for $e$ if $y$ is constant.
\section{Problem Formulation}

Similarly to (\ref{CS}) the optimization problem in (\ref{QCS}) is not convex due to the $\|\mathbf{x}\|_0$ objective function. Therefore first we introduce a convex relaxation to (\ref{QCS}). First note that we can rewrite the (\ref{QCS}) into a general quadratic form

\begin{equation}
y_i = \begin{bmatrix}1 \;\;\; \mathbf{x}^H \end{bmatrix} \begin{bmatrix}a_i \;\;\; \mathbf{b}^H_i \\ \mathbf{c}_i \;\;\; \mathbf{Q}_i \end{bmatrix} \begin{bmatrix}1 \\ \mathbf{x} \end{bmatrix} \in \mathbb{C}, i=1,\ldots,N.
\end{equation}

Since $y_i$ is a scalar using the trace operator we have

\begin{equation}
\begin{split}
y_i =& \text{Tr}(\begin{bmatrix}1 \;\;\; \mathbf{x}^H \end{bmatrix} \begin{bmatrix}a_i \;\;\; \mathbf{b}^H_i \\ \mathbf{c}_i \;\;\; \mathbf{Q}_i \end{bmatrix} \begin{bmatrix}1 \\ \mathbf{x} \end{bmatrix}) \\
=& \text{Tr}(\begin{bmatrix}a_i \;\;\; \mathbf{b}^H_i \\ \mathbf{c}_i \;\;\; \mathbf{Q}_i \end{bmatrix} \begin{bmatrix}1 \\ \mathbf{x} \end{bmatrix} \begin{bmatrix}1 \;\;\; \mathbf{x}^H \end{bmatrix}) \\
=& \text{Tr}(\mathbf{\Phi}_i \mathbf{X})
\end{split}
\end{equation}

where we define $\mathbf{\Phi}_i = \begin{bmatrix}a_i \;\;\; \mathbf{b}^H_i \\ \mathbf{c}_i \;\;\; \mathbf{Q}_i \end{bmatrix}$ and $\mathbf{X} = \begin{bmatrix}1 \\ \mathbf{x} \end{bmatrix} \begin{bmatrix}1 \;\;\; \mathbf{x}^H \end{bmatrix}$. By definition $\mathbf{X}$ is a Hermitian matrix and it satisfies the constraints $\mathbf{X}_{1,1}=1$ and $\text{rank}(\mathbf{X})=1$. The optimization problem (\ref{QCS}) can then be rewritten as

\begin{equation}
\begin{split}
\text{min}_\mathbf{X}\;\; &\|\mathbf{X}\|_0 \\
\text{subject to}\;\; & y_i=\text{Tr}(\mathbf{\Phi}_i\mathbf{X}) , \; i=1,\ldots ,N \\
&\text{rank}(\mathbf{X})=1,\; \mathbf{X}_{1,1}=1,\; \mathbf{X} \succeq 0. \label{QCS2}
\end{split}
\end{equation}

Recasting the optimization problem to solve for the matrix $\mathbf{X}$ in this form is known as matrix lifting \cite{Mcom}, and it was shown that solving for $\mathbf{X}$ one can obtain $\mathbf{x}$ by the rank 1 decomposition of $\mathbf{X}$ by the singular value decomposition (SVD). The problem (\ref{QCS2}) is still not convex, and the convex approximation is given by

\begin{equation}
\begin{split}
\text{min}_\mathbf{X}\;\; &\text{Tr}(\mathbf{X}) + \lambda \|\mathbf{X}\|_1 \\
\text{subject to}\;\; & y_i=\text{Tr}(\mathbf{\Phi}_i\mathbf{X}) , \; i=1,\ldots ,N \\
&\mathbf{X}_{1,1}=1,\; \mathbf{X} \succeq 0. \label{QBP}
\end{split}
\end{equation}

Here, the $\|\mathbf{X}\|_1$ denotes the element-wise $\ell_1$-norm, or the sum of magnitudes of all elements in $\mathbf{X}$, and $\lambda$ is a design parameter. The trace of $\mathbf{X}$ is known to be the convex surrogate of the low rank constraint, and the $\|\mathbf{X}\|_1$ is the convex surrogate of $\|\mathbf{X}\|_0$. The problem (\ref{QBP}) is referred to as quadratic basis pursuit (QBP) and the remainder of this report will be based on solving this optimization problem.

\section{Restricted Isometry Property}

We provide a solution of an undetermined Linear System. consider a linear problem $y = \Phi x$ in which $\Phi$ is a a
$m \times N$ where $m$ is the number of measurement and $m<N$, in
which. Here $x$ will have many solution but we require the
sparsest solution. The formulation is
\begin{equation}
\begin{split}
\text{min}_\mathbf{x}\;\; & \|\mathbf{x}\|_0 \\
\text{subject to}\;\; & \mathbf{y} = \Phi \mathbf{x} \label{CS}
\end{split}
\end{equation}
The solution thus obtained will be the sparse solution for the linear equation $y = \Phi x$ but the above optimization problem will be a NP hard problem.
\newline
\newline
The method described above is very computationally expensive. Eq(1) is a NP-hard problem and can not be solved optimization algorithm hence we can construct other ways of obtaining sparse solution like following problems:
\begin{equation}
\begin{split}
\text{min}_\mathbf{x}\;\; & \|\mathbf{x}\|_1 \\
\text{subject to}\;\; & \mathbf{y} = \Phi \mathbf{x} \label{CS}
\end{split}
\end{equation}
or

\begin{equation}
\begin{split}
\text{min}_\mathbf{x}\;\; & \|\mathbf{x}\|_q \\
\text{subject to}\;\; & \mathbf{y} = \Phi \mathbf{x} \label{CS}
\end{split}
\end{equation}
\subsection{Restricted Isometry Property}

Consider a Linear  equation, $\mathbf{y} = \Phi \mathbf{x}= \Phi \mathbf{x}_{0}$ now using equation (1) to obtain $\mathbf{x}$ from $\mathbf{y}$ that will be equal to or close to $x_0$. The condition for the equality is derived using the restricted Isometry property. Further restricted Isometry property of $\Phi$ prevents any Amplification of  Noise.
A vector is said to be s-sparse if it has at most $s$ nonzero entries. Consider now equation (1), we will explain how the restricted Isometry is applied to (1) to construct a constrain on $\Phi$ in order to maintain the uniqueness of the solution. From (6) we know that $\mathbf{y} = \Phi \mathbf{x}= \Phi \mathbf{x}_{0}$, and  we know that $\mathbf{x}$ can be derived from $\mathbf{y}$ using equation (1), lets  consider  that $x*$  is the  solution for  (1) along with $x_0$.

\section{The $l_1$ approach}
Although the problem in (1) needs a non polynomial time to solve, in general it can be much more effectively solved using $l_1$ minimisation approach. Let us review this approach in this approach in the section. The $l_1$ minimization problem is the following,
\begin{equation}
\begin{split}
\text{min}_\mathbf{x}\;\; & \|\mathbf{x}\|_1 \\
\text{subject to}\;\; & \mathbf{y} = \Phi \mathbf{x} \label{CS}
\end{split}
\end{equation}

Since  the $l_1$ minimization is a convex problem and can be converted easily into linear programming problem it can be easily  solved.  Now we have to see under what condition does the solution of (1) is same as the solution of (11). we can define the condition by too Concept: 1) Mutual Coherence, 2) Restricted Isometric property. In this paper we will provide solution using Restricted Isometric property.

\subsection{The Exact recovery condition for $l_1$ norm}
In this section some, theoretical results are presented on uniqueness and recoverability of QBP in the noiseless case. Proofs are omitted in this report but can be found in \cite{main}. For convenience we first introduce a linear operator $B:$

\begin{equation}
B:\; X\in \mathbb{C}^{n\times n} \mapsto {\text{Tr}(\mathbf{\Phi}\mathbf{X})}_{1\leq i\leq N} \in \mathbb{C}^{N} \label{B}
\end{equation}
using (\ref{B}) we can define the restricted isometry property.

\newtheorem{Definition}{Definition}
\newtheorem{Theorem}{Theorem}
\begin{Definition}[\textbf{RIP}]
A linear operator $B(\cdot)$ is ($\epsilon,k$)-RIP if

\begin{equation}
\left| \frac{\|B(\mathbf{X})\|^2_2}{\|\mathbf{X}\|^2_2} \right| < \epsilon
\end{equation}
for all $\|\mathbf{X}\|_0\leq k$ and $\mathbf{X}\neq0$.
\end{Definition}

From \cite{main} it we have the following theorem for uniqueness,

\begin{Theorem}[\textbf{Uniqueness}]
$\overline{\mathbf{x}}$ is a solution to (\ref{QCS}) and if $\mathbf{X}^*$ satisfies $y=B(\mathbf{X}^*)$, $\mathbf{X}^* \succeq 0$, $\text{rank}(\mathbf{X}^*)=1$, $\mathbf{X}^*_{1,1}=1$ and if $B(\cdot)$ is a $(\epsilon,2\|\mathbf{X}^*\|_0)$-RIP with $\epsilon<1$ then the rank-1 decomposition of $\mathbf{X}^*$ is equal to $[1\;\;\; \overline{x}]^T$.
\end{Theorem}

The next concern is on recoverability of the QBP algorithm, first we define the mutual coherence of a dictionary matrix $\mathbf{A}$

\begin{Definition}[Mutual Coherence]
For a matrix $\mathbf{A}$ the mutual coherence $\mu(\mathbf{A})$ is given by

\begin{equation}
\mu(\mathbf{A})=\text{max}_{i,j} \frac{|\mathbf{A}^H_i \mathbf{A}_j|}{\|\mathbf{A}_i\|_2\|\mathbf{A}_j\|_2} , i\neq j.
\end{equation}
\end{Definition}
Also define a matrix $\mathbf{D}$ be a matrix satisfying $\mathbf{y}=\mathbf{D}\mathbf{X}^s=B(\mathbf{X})$ where $\mathbf{X}^s=\text{vec}(\mathbf{X})$. Then the rank-1 decomposition of the solution of QBP is equal to $[1\;\;\; \overline{\mathbf{x}}]^T$ if $\|\mathbf{X}\|_0<0.5(1+1/\mu(D)).$

\subsection{Stability}

In the previous section we proved condition for exact
recovery for $l_1$ norm now we move further and prove the
condition for stability. When we apply an undetermined Linear
transform $H$ to a signal $x_0$ there could be some error induced
in the undetermined linear system or x0 could have noise in
it . In such cases we need to guarantee that the x calculated
from measurement y will not blow up in comparison to x0
For this there are two set of theorem defining the upper limit
of error between original signal and the recovered signal. we
consider that y is contaminated and incomplete observation
$y = Hx_0 + e$, e is the error term .

\section{The $l_q$ approach}

In this section some, theoretical results are presented on uniqueness and recoverability of QBP in the noiseless case. Proofs are omitted in this report but can be found in \cite{main}. For convenience we first introduce a linear operator $B:$

\begin{equation}
B:\; X\in \mathbb{C}^{n\times n} \mapsto {\text{Tr}(\mathbf{\Phi}\mathbf{X})}_{1\leq i\leq N} \in \mathbb{C}^{N} \label{B}
\end{equation}
using (\ref{B}) we can define the restricted isometry property.

\begin{equation}
\left| \frac{\|B(\mathbf{X})\|^2_2}{\|\mathbf{X}\|^2_2} \right| < \epsilon
\end{equation}
for all $\|\mathbf{X}\|_0\leq k$ and $\mathbf{X}\neq0$.

A linear operator $B(\cdot)$ is ($\epsilon,k$)-RIP if

\begin{equation}
\left| \frac{\|B(\mathbf{X})\|^2_2}{\|\mathbf{X}\|^2_2} \right| < \epsilon
\end{equation}
for all $\|\mathbf{X}\|_0\leq k$ and $\mathbf{X}\neq0$.

From \cite{main} it we have the following theorem for uniqueness,

\begin{Theorem}[\textbf{Uniqueness}]
$\overline{\mathbf{x}}$ is a solution to (\ref{QCS}) and if $\mathbf{X}^*$ satisfies $y=B(\mathbf{X}^*)$, $\mathbf{X}^* \succeq 0$, $\text{rank}(\mathbf{X}^*)=1$, $\mathbf{X}^*_{1,1}=1$ and if $B(\cdot)$ is a $(\epsilon,2\|\mathbf{X}^*\|_0)$-RIP with $\epsilon<1$ then the rank-1 decomposition of $\mathbf{X}^*$ is equal to $[1\;\;\; \overline{x}]^T$.
\end{Theorem}

The next concern is on recoverability of the QBP algorithm, first we define the mutual coherence of a dictionary matrix $\mathbf{A}$

\begin{Definition}[Mutual Coherence]
For a matrix $\mathbf{A}$ the mutual coherence $\mu(\mathbf{A})$ is given by

\begin{equation}
\mu(\mathbf{A})=\text{max}_{i,j} \frac{|\mathbf{A}^H_i \mathbf{A}_j|}{\|\mathbf{A}_i\|_2\|\mathbf{A}_j\|_2} , i\neq j.
\end{equation}
\end{Definition}
Also define a matrix $\mathbf{D}$ be a matrix satisfying $\mathbf{y}=\mathbf{D}\mathbf{X}^s=B(\mathbf{X})$ where $\mathbf{X}^s=\text{vec}(\mathbf{X})$. Then the rank-1 decomposition of the solution of QBP is equal to $[1\;\;\; \overline{\mathbf{x}}]^T$ if $\|\mathbf{X}\|_0<0.5(1+1/\mu(D)).$

\subsection{Stability}
In the previous section we proved condition for exact
recovery for $l_1$ norm now we move further and prove the
condition for stability. When we apply an undetermined Linear
transform $H$ to a signal $x_0$ there could be some error induced
in the undetermined linear system or x0 could have noise in
it . In such cases we need to guarantee that the x calculated
from measurement y will not blow up in comparison to x0
For this there are two set of theorem defining the upper limit
of error between original signal and the recovered signal. we
consider that y is contaminated and incomplete observation
$y = Hx_0 + e$, e is the error term .

\section{Algorithm}

In this section we demonstrate the use of QBP. We draw $\mathbf{B} \in \mathbb{C}^{N\times n}$, and $\mathbf{Q}_i \in \mathbb{C}^{n\times n},\; i= 1,\ldots N$ from a complex Gaussian random variable. We set $\mathbf{C}=0$ and $\mathbf{a}=0$, $N=18$, $n=20$, and $\mathbf{x}$ is real and shown in Fig. 1. All scenarios are the noiseless case.

\begin{figure}
\centering
\includegraphics[width=1\columnwidth]{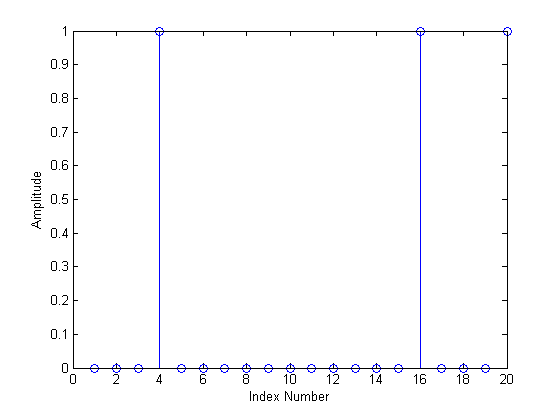}
\caption{The ground truth}
\label{fig:truth}
\end{figure}

Fig. 2 shows the output of using the basis pursuit de-noising algorithm which takes a similar form as (\ref{BP}). The optimization problem that was solved instead of (\ref{BP}) was,

\begin{equation*}
\text{min}\;\; \|\mathbf{y-Bx}\|^2_2 + \lambda \|\mathbf{x}\|_1.
\end{equation*}

This optimization problem was solved rather than (\ref{BP}) because it was found that solving (\ref{BP}) often led to the result that the problem is infeasible. The regularization parameter $\lambda$ was chosen to be 50. As expected since basis pursuit algorithms can only account for the linear portion of the model it performs badly and the estimate is far from the truth.

\begin{figure}
\centering
\includegraphics[width=1\columnwidth]{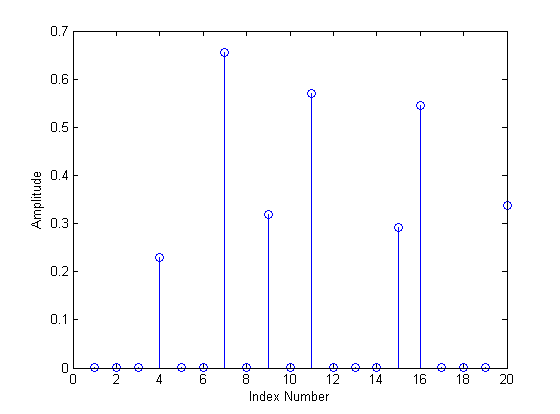}
\caption{Output of the basis pursuit de-noising  algorithm with $\lambda = 50$, estimate is far from the ground truth}
\label{fig:CS}
\end{figure}

Fig. 3 shows the result of the QBP with $\lambda=50$, the result is a perfect reconstruction and QBP can account for both the quadratic and the linear terms.

\begin{figure}
\centering
\includegraphics[width=1\columnwidth]{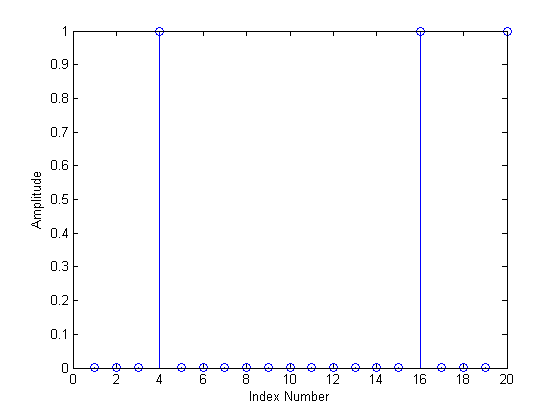}
\caption{Output of the QBP algorithm with $\lambda=50$, perfect reconstruction is obtained}
\label{fig:Q_CS50}
\end{figure}

Fig. 4 shows the result of the QBP without applying a regularization parameter or setting $\lambda=0$, this is the same as not enforcing the sparsity constraint. It is seen that without the sparsity constraint we recover a very dense estimate that and many solutions can exist.

\begin{figure}
\centering
\includegraphics[width=0.5\columnwidth]{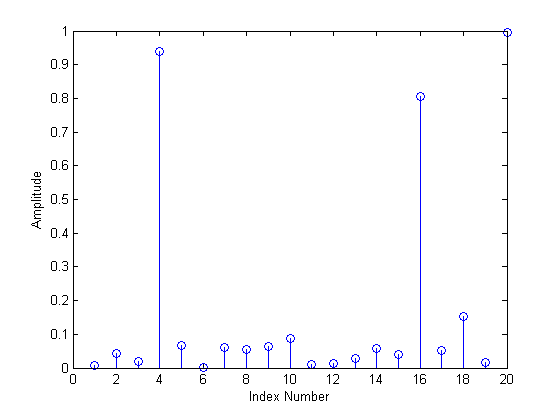}
\caption{Output of the QBP algorithm with $\lambda=0$, dense estimate is obtained}
\label{fig:Q_CS0}
\end{figure}

\section{Results}
\begin{figure}
\centering
\includegraphics[width=0.9\columnwidth]{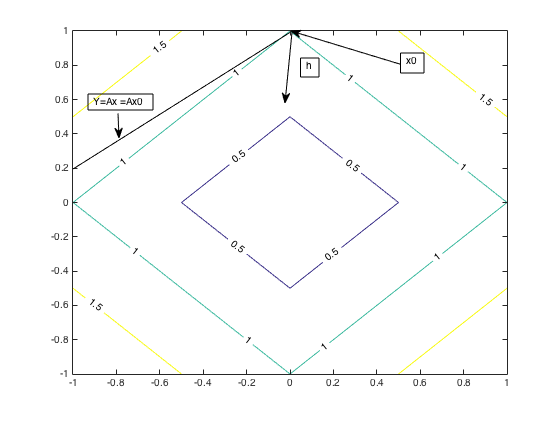}
\caption{Describing the geometric shape of y =Ax = Ax0 for successful
recovery of x from y such that $x=x_0$}
\label{fig:CS}
\end{figure}

\begin{figure}
\centering
\includegraphics[width=0.9\columnwidth]{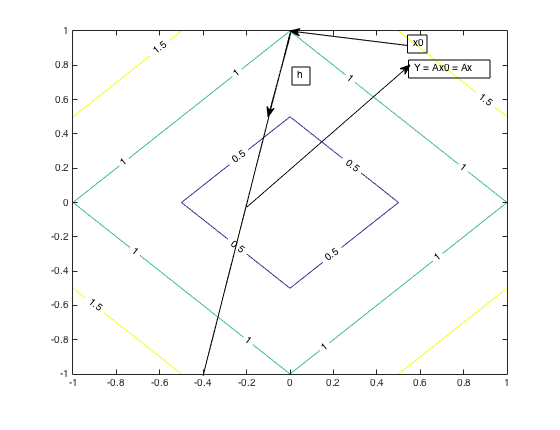}
\caption{Describing the geometric shape of $y =Ax = Ax_0$ for recovery of x
from y such that $x \neq x_0$ and hence no exact recovery}
\label{fig:CS}
\end{figure}

\begin{figure}
\centering
\includegraphics[width=0.9\columnwidth]{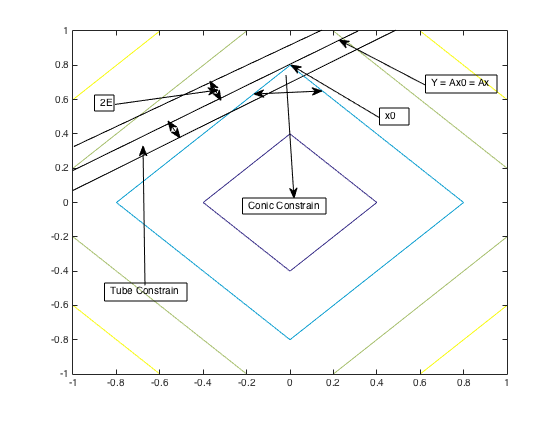}
\caption{Describing the geometric shape of $Y = Ax+e$ and the Tube constrain,
and conic constrain induced due to $l_1$ minimization}
\label{fig:CS}
\end{figure}

\begin{figure}
\centering
\includegraphics[width=0.9\columnwidth]{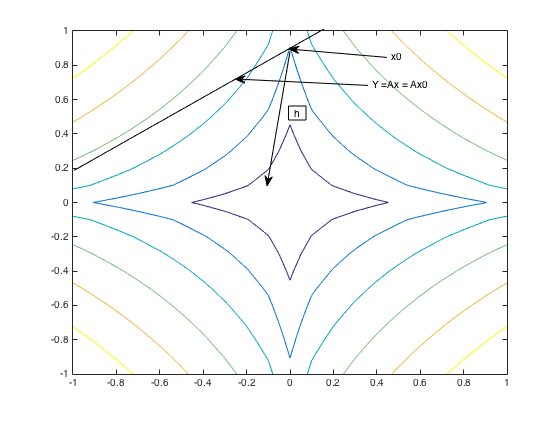}
\caption{Describing the geometric shape of $Y = Ax = Ax_0$ for successful
recovery of x from y such that $x = x_0$ . The decoder used here is $\delta_q$ and x is
minimized over optimisation problem}
\label{fig:CS}
\end{figure}

\begin{figure}
\centering
\includegraphics[width=0.9\columnwidth]{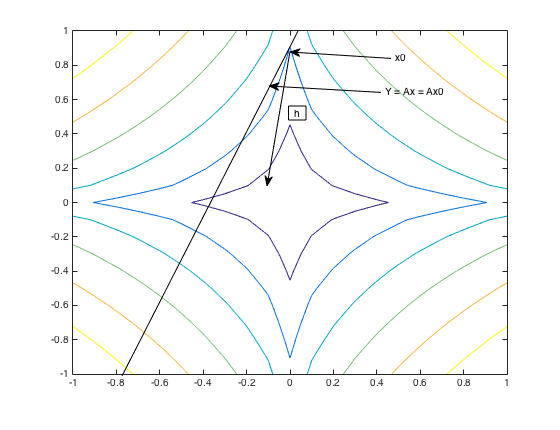}
\caption{Describing the geometric shape of $Y = Ax = Ax_0$ for successful
recovery of x from y such that $x = x_0$ . The decoder used here is $\delta_q$ and x is
minimized over optimisation problem}
\label{fig:CS}
\end{figure}

\begin{figure}
\centering
\includegraphics[width=0.9\columnwidth]{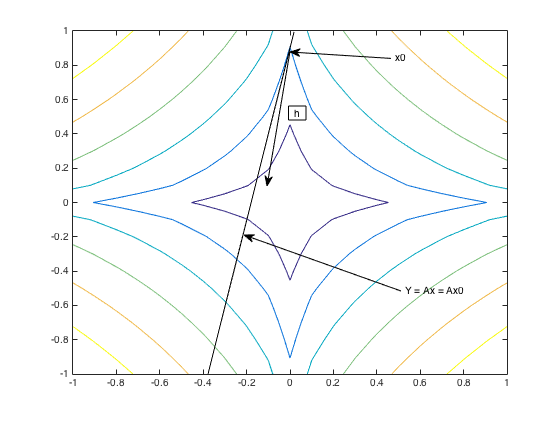}
\caption{Describing the geometric shape of $Y = Ax = Ax_0$ for successful
recovery of x from y such that $x = x_0$ . The decoder used here is $\delta_q$ and x is
minimized over optimisation problem}
\label{fig:CS}
\end{figure}

\begin{figure}
\centering
\includegraphics[width=0.9\columnwidth]{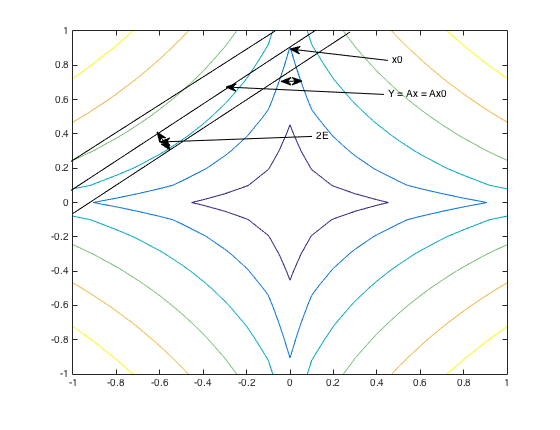}
\caption{Describing the geometric shape of $Y = Ax = Ax_0$ for successful
recovery of x from y such that $x = x_0$ . The decoder used here is $\delta_q$ and x is
minimized over optimisation problem}
\label{fig:CS}
\end{figure}
\section{Conclusion}
This paper helped understand that geometrically
description of a basic Pursuit problem and I can use a
similar theory for making a $l_q$ minimisation problem have a
global minima.


\begin{thebibliography}{1}
\bibitem{Mcom4} Restricted Isometry and its implication for compresses sensing by Emmanuel J. Candes
\end{thebibliography}
\end{document}